\documentclass[conference,10pt]{IEEEtran}

\usepackage[T1]{fontenc}
\usepackage[dvips]{graphicx}
\usepackage[english]{babel}
\usepackage[latin1]{inputenc}

\usepackage{amssymb}
\usepackage{amstext}
\usepackage{url}

\usepackage{cite} 
\usepackage{stfloats}  
                        
\begin{document}

\title{Time Allocation of a Set of Radars in a Multitarget Environment}
\author{ \authorblockN{Emmanuel Duflos}
  \authorblockA{
    LAGIS UMR CNRS 8146\\
    INRIA Futurs\\
    Ecole Centrale de Lille
    BP 46\\
    59651 LILLE Cedex\\
    FRANCE\\
    Email: emmanuel.duflos@ec-lille.fr} \and \authorblockN{Marie de
    Vilmorin} \authorblockA{
    FRANCTechno parc futura\\
    Universi\'e d'Artois\\
    Facult\'e des Sciences Appliqu\'ees\\
    62400 Bethune\\
    FRANCE\\
    Email: marie.de\_vilmorin@fsa.univ-artois.fr} \and
  \authorblockN{Philippe Vanheeghe}
  \authorblockA{LAGIS UMR CNRS 8146\\
    INRIA Futurs\\
    Ecole Centrale de Lille
    BP 46\\
    59651 LILLE Cedex\\
    FRANCE\\
    Email: philippe.vanheeghe@ec-lille.fr} }

\maketitle

\selectlanguage{english}

\begin{abstract}
The question tackled here is the time allocation of radars in a multitarget environment. At a given time radars can only observe a limited part of the space; it is therefore necessary to move their axis with respect to time, in order to be able to explore the overall space facing them. Such sensors are used to detect, to locate and to identify targets which are in their surrounding aerial space. In this paper we focus on the detection schema when several targets need to be detected by a set of delocalized radars. This work is based on the modelling of the radar detection performances in terms of probability of detection and on the optimization of a criterion based on detection probabilities. This optimization leads to the derivation of allocation strategies and is made for several contexts and several hypotheses about the targets locations.

\end{abstract}


%
\section{Introduction}
\label{sec:Introduction}

In many applications sensors are nowadays a part of a multisensor
system, each sensor bringing its complementarity and its redundancy to
the overall system. Year after year the complexity and the
performances of many sensors have increased leading to more and more
complex multisensor systems which supply the decision centers with an
increasing amount of data. This increasing complexity also led to
other uses for each sensor and therefore for the multisensor systems.
It is no more considered as a passive system the role of which is just
limited to simple measurement actions; the many parameters of each
sensor and the interactions between all the sensors allow to choose how
the measurement action must be done: the sensors need to be managed.
The complexity of this problem is such that it is often impossible to
a man to find an optimal solution (with respect to the goal of the
mission of the multisensor system) and multisensor management
strategies must be derived. That is the reason why sensors management
has become during the past years an active field of research. From a
theorical point of view this problem can be written in the frame of
optimal control and the sensor management viewed as a Markov decision
problem. Optimal solutions could therefore be found. Unfortunately,
the complexity is such that it is impossible in practice to derive
these solutions. Sub-optimal solutions as well as alternative
approaches have then been proposed. In \cite{kreucher2004} or
\cite{kreucher2005} the authors use reinforcement learning, Q-learning
and approximation functions to derive sub-optimal solutions. In many
works the choice of the next action is based on information theory and
information divergence like the R\'{e}nyi information divergence and the
Kullback Leibler divergence \cite{kastella1997},
\cite{kreucher2004},\cite{kreucher2005}. In \cite{mahler1996}
Mahler proposes to solve the problem in the frame of random sets. All
these works bring a possible solution to the sensor management problem
but as far as the authors know, it is often difficult to derive bound
of performance which can be a drawback in an operational context.
Moreover, these approaches rarely take into account the
characteristics of the sensors. The work described in this paper
proposes, in the frame of an aerial patrol in charge of the detection
of potential targets, to derive radar optimal time allocations (a part
of the sensor management problem) which allow to determine such bounds
and which are based on the modelling of the detection performances of
a radar. It is assumed here that each aircraft is equipped with an ESA
(Electronically Steered Antenna) radar. We focus on the detection
step for which a fixed duration $T$ has been allocated. Methods exist
to optimize the detection of a single target by a single sensor 
and the frame Search Theory is devoted to such a
problem \cite{BLA99}, \cite{lecadre2000}. In this paper we consider a multitarget
environment and the optimization process is led by considering the
overall targets and not the targets one by one. The problem then
becomes: if radars have to observe $P$ targets during $T$, how do they
organize themselves to detect them in the best possible way, i.e. how
do they distribute the duration $T$ over the space directions ? The
aim of this article is to derive an optimal temporal allocation based
on the modelling of the radar detection probability and on an \textit{a
  priori} knowledge coming from an ESM type (Electrical Support
Measurement) or AEW type (Airbone Early Warning) system of supervision. 
Along the study, two contexts are considered. The first one is
the ideal case: the position of the targets are known and we must
detect them. Of course this situation is not realistic but it allows
to derive some interesting results for the second context : the
position of the targets are known by the mean of probability
densities. After having defined the assumptions of our study in the
second section, we present in the third section a modelling of the
radar detection functions. However the context of this study is
multisensor multitarget, we start by a study of the optimization of
the detection process in a monosensor monotarget environment.
Comparing to existing methods, our aim in this preliminary work is to
derive analytically an optimal strategy and the corresponding
probability of detection. This last probability will be used along the
overall paper. The third section presents analytic results and a
performance evaluation. The multitarget environment is tackled in the
fourth section but we are still in a monosensor case. Under the
assumption of an \textit{a priori} knowledge, we propose an optimal
temporal allocation. The allocation derived in this section uses the
results derived in the previous sections. Finally, the last section
shows how all the previous results can be used to propose an
allocation strategy in the multisensor multitarget case. It is
important to understand the needs at the origin of the study, proposed
by Thales Optronics, the results of which are written out in this
paper. The aim was to found bounds of performance for optimal
allocation strategies. Therefore, the fact of considering
deterministic knowledge about the targets, as it is the case in some
parts of the paper, has sense even if it is not realistic in an real
operational context. When used or adapted in such a context, the
proposed strategies are no more optimal but we know the bounds of
performance which can be interesting.

\section{Hypotheses}

The main assumption of this article is the use of an \textit{a priori}
knowledge. By this expression we mean a knowledge about the situation,
in particular a knowledge about the positions of the targets. This
assumption is justified by the integration of the sensors in a
supervision system of the ESM type (Electrical Support Measurement) or
AEW type (Airbone Early Warning) for instance. Using these sensors, it
is possible to obtain information on the angular positions of the
targets and then to derive information on their distances from the
sensors. In this paper, the \textit{a priori} knowledge is ideal or
deterministic (for reasonning purposes) or more realistic (given by
density probability function). We also consider a 2D space; this
assumption does not limit the general character of the study, however
it reduces the calculations. Finally, we suppose that the observation
durations are sufficiently short so that the aircrafts can be regarded
as stationary which means that the targets do not move out their
resolution cell during the observation process.

\section{Detection probability optimization in a monosensor monotarget
  context
\label{capteur}}

\subsection{The radar sensor}
\label{radar}

In order to establish optimal management strategies for the sensors,
it is necessary to understand their operating mode. In particular, the
modelling of the detection probability is a fundamental basis for the
management strategies we are going to define. The radar is an active
sensor since it emits a signal which is reflected on the target. The
radar considered here has an electrical scanning. It means that its
mecanical axis is fixed and that it is the direction of the analyzing
wave which is modified during the observation. First, we are interested 
in the signal to noise ratio. If the sensor observes during
$T$ a target at the range $r$ in a direction which forms an angle
$\theta $ with the mechanical axis of the antenna, then the signal to
noise ratio of an echo is equal to \cite{DI04}:
\begin{equation}
  SNR=\frac{\alpha T\cos ^{2}\left( \theta \right) }{r^{4}}  \label{eq1}
\end{equation}
where $\alpha $ is an operational parameter, which depends of the
radar and the target (target's radar cross section). In this paper,
the targets are supposed to be similar, i.e. to have the same
reflexion power, then $\alpha $ is constant. This expression is
established in a context where the disturbing signal, which is
supposed to be only due to the thermal noise of the radar, is modelled
by a normal random variable. Once the expression of the signal to
noise ratio known, the target detection probability can be derived. It
is modelled as:

\begin{equation}
  P_{d}=\left( P_{fa}\right) ^{\frac{1}{1+SNR}}  \label{eq2}
\end{equation}

with $P_{fa}$ the false alarm probability which is taken here to one
false alarm per second and per resolution cell. This expression is
established under the assumptions of a fluctuating target and a
modelling of the received energy of type ''Swerling 1'' \cite{KLE97}.
A target does not have a regular form, therefore the reflected energy
varies from an impulse to another. The target can then be considered
as a set of elementary reflectors the positions of which in the space
are related to the target orientation. The returned signals are then
independent and the amplitude of the received energy fluctuates. These
targets are called ''fluctuating targets'' and have been modelled by
Swerling \cite{SWE60} : they are called \textit{Swerling 0, Swerling
  1, Swerling 2, Swerling 3 } and \textit{Swerling 4}. The
\textit{Swerling 1 } type is particularly adapted to the case of the
air target detection.

\subsection{Optimization of the detection probability in a monosensor
  monotarget context}

According to expressions (\ref{eq1}) and (\ref{eq2}) it can be easily
shown that the detection probability is strongly degraded when the
range increases. There are several methods to improve it.  A first
solution is to use a procedure of "alert and confirmation"
\cite{DAN81},\cite{BLA99}. This method consists in doing two detection
steps: the first one with a low detection threshold, the second one
with a higher one in order to eliminate false alarms from the alert
step. During this second step, the emitted wave is adapted to the
target. However, this process of decomposition in two steps needs a
long integration time. A solution could be to increase it but for high
detection probability, the slope $\frac{dP_{d}}{dt}$ is small. Instead
of carrying out only one acquisition of the signal during $T,$
therefore only one detection, we propose to acquires $N$ elementary
signals and to carry out an elementary detection on each received
signal, that is to realize $N$ elementary detections. The use of the
radar with a different emission frequency at each elementary
detection allows to obtain independent detections and allows the
analytic derivation of an optimal detection probability as it is shown
in the following \cite{DI04}, \cite{STI98}. If $P_{de}$ denotes the
elementary detection probability, then the cumulative detection
probability is equal to \cite{KLE97}:
\begin{equation}
  P_{d}=1-\left( 1-P_{de}\right) ^{N}  \label{eq3}
\end{equation}
where $P_{de}$ is given by the expression (\ref{eq2}), for an observation
duration equal to $\frac{T}{N}.$ The problem is to find the number $N$ of elementary detections which optimizes this
cumulative detection probability. By considering the target's signal to noise ratio far higher than one,
it is possible to detail the expression of the elementary detection
probability as:
\begin{equation}
  P_{de}=\left( P_{fa}\right) ^{\frac{1}{SNR}}=\exp \left( \frac{r^{4}N\ln
      \left( P_{fa}\right) }{\alpha T\cos ^{2}\left( \theta \right) }\right)
  \label{eq4}
\end{equation}

At this point it is important to understand what is the meaning and the limitation of "the signal to noise ratio is far higher than one". It means that it is high enough to make the approximation of $P_{de}$ by (\ref{eq4}), but it is not high enough to consider that this probability is almost equal to one; a detection phase is therefore necessary. We now introduce the constant $\beta$ as :

\begin{equation}
\beta =\frac{r^{4}\ln \left(\frac{1}{P_{fa}}\right) }{\alpha T\cos ^{2}\left( \theta \right)}
\end{equation}

which allows us to express the probability (\ref{eq4}) as $P_{de}=\exp \left( -\beta N\right)$. The cumulative detection probability is then equal to:
 
\begin{equation}
P_{d}=1-\exp \left( N\ln \left( 1-\exp \left( -\beta N\right) \right) \right).
\label{eq5}
\end{equation}

Using a classical optimization process on this last probability, it can be shown that it is optimal if $N$ is equal to :

\begin{equation}
  N_{opt}=\frac{\gamma _{r}\alpha T\left( \cos \left( \theta \right) \right)
    ^{2}}{r^{4}\ln \left( \frac{1}{P_{fa}}\right) }  
\label{eq7}
\end{equation}

with $\gamma _{r}=-\beta N=-\ln \left( P_{de}\right)=\ln 2$. The elementary detection probability $P_{de}$ is therefore equal to $0.5$ and the cumulative one is equal to:

\begin{equation}
  P_{d}=1-\exp \left( -\frac{T}{\tau _{r}}\right)  \label{eq8}
\end{equation}

with

\begin{equation}
  \tau _{r}=\frac{r^{4}\ln \left( P_{fa}\right) }{\gamma _{r}\alpha \left(
      \cos \left( \theta \right) \right) ^{2}\ln \left( 1-\exp \left( -\gamma
        _{r}\right) \right) }.  \label{eq8a}
\end{equation}

These results show on the one hand that the modelling of the radar
sensor detection functions makes it possible the elaboration of
analytical strategies of optimization of the detection probability
and, on the other hand, that it is possible to quantify the
performances. A few remarks about these results:

\begin{itemize}
\item The optimal number of elementary detections is not a natural, $%
  N_{opt}\notin \mathbb{N}$. However, it does not alter the general frame
  of our method and allows us to calculate optimal performances which
  will be used as references, like the Cram\'er-Rao lower bound in
  estimation theory.

\item The assumption of an important signal to noise ratio is a trick
  which allows us to write the detection probability simplier.
  However, the probabilities obtained can be close to $0.5,$ which
  justifies the elaboration of an optimal detection process.
\end{itemize}

In the next section we will see how to use these results to detect several targets.

\section{Monosensor Multitarget Environment}
\subsection{Deterministic knowledge
\label{hyp1}}

We consider a situation where $P$ targets are present in an air space.
The knowledge about them is such that their angular deviations $\theta
_{i}$ and their ranges $r_{i}$ are known $\forall $ $%
i\in \left\{ 1,..,P\right\}$. Our goal is to detect them with a radar.
Furthermore, targets are supposed to be localized in different
directions of the space, that is $\theta _{i}\neq \theta _{j}$
$\forall \left( i,j\right) \in \left\{ 1,..,P\right\} ^{2}.$ Finally,
we also suppose that all the targets don't represent the same threat
which leads to the introduction of weights $\varepsilon_i$. In these
conditions, let us call $t_{i}$ the observation duration of the target $i$
by the radar. Since angular deviations are different, the radar
does not observe several targets simutalneously and it results in the
following relation between durations $t_{i}$ and the total duration
$T$:
\begin{equation}
\sum_{i=1}^{P}t_{i}=T \text{ with } t_{i}\geq 0\quad \forall i\in \left\{ 1,..,P\right\}.
\label{eq9}
\end{equation}
Durations must obviously be positive, but we will see that they could
be null, because of the relative positions of the targets with respect
to the sensor. Our aim is to maximize the detection of all the
targets. As a probability is always positive, maximize each of them is
equivalent to maximize their sum. Then we define the criterion:
 \begin{equation}
 J=\sum_{i=1}^{P}\varepsilon _{i}P_{di}\left( t_{i}\right)  \label{eq11}
 \end{equation}
 where $\varepsilon _{i}$ can be interpreted as a potential threat or
 priority coefficient. This concept of threat is introduced here in a
 general way and its characterization is behind the scope of the
 present paper. These coefficients can for instance be inversely
 proportionnal to the distance \cite{VAN01}. $P_{di}\left( t_{i}\right)$ 
 is the detection probability of the target $i$ by the radar, for a
 duration $t_{i}$. If this duration is known, then we are in the context
 described in section \ref{capteur}: a monosensor monotarget context.
 Using previous results it is possible to write:
\begin{equation}
P_{di}\left( t_{i}\right) =1-\exp \left( -\frac{t_{i}}{\tau _{r_{i}}}\right)
\label{eq12}
\end{equation}
with 
\begin{equation}
\tau _{ri}=\frac{r^{4}\ln \left( P_{fa}\right) }{\gamma _{r}\alpha \left(
\cos \left( \theta _{i}\right) \right) ^{2}\ln \left( 1-\exp \left( -\gamma
_{r}\right) \right) }.  \label{eq13}
\end{equation}
 According to this last expression of the detection probability, the criterion $J$
 will reach its maximum when the durations $t_{i}$ tend towards infinite,
 which is not compatible with the temporal constraint which was define. Our
 aim is then to optimize the criterion (\ref{eq11}) under the constraints (\ref{eq9}): 
 \begin{equation}
 \mathcal{J}\left\{ 
 \begin{array}{ll}
 \text{maximize} & 
 \begin{array}{c}
 J=\sum\limits_{i=1}^{P}\varepsilon _{i}P_{di}\left( t_{i}\right)
 \end{array}
 \\ 
 \text{under the constraints} & \left\{ 
 \begin{array}{l}
 \sum\limits_{i=1}^{P}t_{i}=T\quad \\ 
 t_{i}\geq 0\quad \forall i\in \left\{ 1,..,P\right\}
 \end{array}
 \right.
 \end{array}
 \right.  \label{eq14}
 \end{equation}

 With such a formulation, we face to a classical problem of
 optimization under constraints which can be solved by the way of
 Lagrangian functions and tools of convex optimization. Such an
 optimization process leads to the following results. Let us introduce
 the function $x\mapsto \left\lfloor x\right\rfloor ^{+}$ defined on
 $\mathbb{R}$ by:

 \begin{equation}
 \begin{array}{lll}
 \left\lfloor x\right\rfloor ^{+} & = & x\;if\;x>0 \\ 
 & = & 0\;else
 \end{array}
 \label{eq21}
 \end{equation}

 and $\lambda $ the single solution of the following equation:
 
 \begin{equation}
 \sum_{i=1}^{P}\frac{\tau _{ri}}{T}\left\lfloor \ln \left( \frac{T\varepsilon
 _{i}}{\tau _{ri}\lambda }\right) \right\rfloor ^{+}-1=0.  
 \label{eq22}
 \end{equation}

 If $\mathcal{I}$ is the suffix set defined by:
 
 \begin{equation}
 \mathcal{I}=\left\{ i\in \left\{ 1,..,P\right\} \left| \lambda <%
 \frac{T\varepsilon _{i}}{\tau _{ri}}\right. \right\}   \label{eq23}
 \end{equation}

 then the optimal temporal allocation is given by :
 
 \begin{equation}
 \left\{ 
 \begin{array}{lll}
 t_{i} & = & \tau _{ri}\ln \left( \frac{T\varepsilon _{i}}{\tau _{ri}\lambda }%
 \right) \;if\;i\in \mathcal{I} \\ 
 & = & 0\;else
 \end{array}
 \right.   
 \label{eq24}
 \end{equation}

 Furthermore, the optimal elementary detection number for each $t_{i}$, i.e
 for each target, is equal to: 
 \begin{equation}
 n_{i}=\frac{\gamma _{r}\alpha t_{i}\left( \cos \left( \theta \right) \right)
 ^{2}}{r_{i}^{4}\ln \left( P_{fa}\right) }.  \label{eq25}
 \end{equation}

 Results given in (\ref{eq24}) depend on the parameter $\lambda $ which
 is solution of the equation (\ref{eq22}). This equation has an
 analytical solution if $card\left( \mathcal{I}\right) =P$. In this
 case it can be shown from the previous result that the optimal
 allocation of the duration $T$ between the $P$ targets is given by :

 \begin{equation}
 t_{i}=\frac{\sum_{j=1}^{P}\tau _{rj}\ln \left( \frac{\epsilon _{i}\tau _{rj}%
 }{\epsilon _{j}\tau _{ri}}\right) +T}{\sum_{j=1}^{P}\frac{\tau _{rj}}{\tau
 _{ri}}}\;\forall i\in \left\{ 1,...,P\right\}.  \label{eq39}
 \end{equation}
 A sufficient condition to obtain this result is that there exists $\lambda
 _{0}>\lambda $ such that: 
 \begin{equation}
 \lambda _{0}<\frac{T\epsilon _{i}}{\tau _{ri}}\;\forall i\in \left\{
 1,...,P\right\}  \label{eq40}
 \end{equation}

 So far, we have derived optimal allocations for the detection of $P$
 targets given a total duration. We can remark that an implicite
 assumption has been made: the infinite divisibility of the duration
 $T$, which is not the case in reality. However, this assumption is
 justified in this article by the use of a radar. This sensor having
 an electronic scanning mode, the movement from an angular position to
 another can be considered as instantaneous. The major hypothesis of
 this section is the deterministic knowledge about the situation. We
 propose in the next section an optimal temporal allocation in the
 case of an \textit{a priori} knowledge defined by probability
 densities which corresponds to more realistic context.

 \subsection{Probabilistic knowledge
 \label{sec:probknowledge}}

\begin{figure}
  \centering
  \includegraphics[width=6cm]{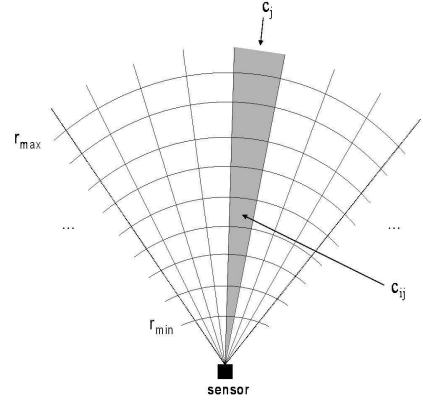}
  \caption{\textit{Representation of the directions and cells of the space}}
  \label{figespace}    
 \end{figure}

 In this section we assume that we have a weak knowledge of the targets
 positions. Several targets can appear in the same direction which is
 more realistic than in the previous section. Thus we have to consider
 all the space and not a few directions as it was possible previously.
 The space directions we consider are the angular fields of view of the
 radar. Moreover, since the sensor observes the globality of a
 direction at the same time, we want to determine the observation
 duration in this direction, that is the duration $t_{j}$ in the
 direction $j$. In order to do so, the observation space is sampled
 according to the resolution cells of the radar as it is described in
 the figure \ref{figespace}.
 The sampled space is therefore defined by the interval of ranges
 $\left[ r_{\min },r_{\max }\right]$ in which the detection is realized
 and the angular sector $c_{j}$ i.e the direction $j,\,j\in \left\{
 1,..,N_{d}\right\}$. We have seen in the section \ref{radar} that
 the sensor forms a set of range resolution cells in each direction,
 i.e in each cell $c_{j}.$ Then we consider a set of sub-cells, the
 resolution cells, $c_{ij},$ at the range $r_{i},$ in the direction
 $j,$ $i\in \left\{1,...,N_{r}\right\}$. Since a sensor observes
 simultaneously all the targets present in the same direction, we are
 going to determine the probability of detecting from one to several
 targets in each of these directions. First we consider the expression
 of the detection probability given by relation (\ref{eq4}). It
 represents the detection probability of a target knowing that it is at
 a range $r$ from the sensor. Let $H_{k}$ be the event "the target $k$ is
 detected". Using the formalism of conditionnal probabilities we can
 write the probability of detecting the target $k$ at a range $r$ as:
 \begin{equation}
 P\left( H_{k},r\right) =P\left( H_{k}\left| r\right. \right) P\left(
 k,r\right)  \label{eq47}
 \end{equation}
 with $P\left( k,r\right) $ the target location probability at the
 range $r.$ We assume that the target detection probability in a given
 cell can be approximated by the detection probability of this target
 at the range on which the cell is centered. The expression of the
 previous probability can therefore be written as:

 \begin{equation}
 P\left( H_{k},c_{ij}\right) =P\left( H_{k}\left| c_{ij}\right. \right)
 P\left( k,c_{ij}\right)  \label{eq48}
 \end{equation}

 where $P\left( k,c_{ij}\right) $ is obtained by the integration of the
 a priori density of probability in the cell $c_{ij}.$ We note $P\left(
   k,c_{ij}\right) =\rho _{ijk}.\,P\left( H_{k}\left| c_{ij}\right.
 \right) $ is derived from the expression (\ref{eq4}), with the target
 at the approximate range $r_{i}$. Finally we obtain:

 \begin{equation}
 P\left( H_{k},c_{ij}\right) =P_{dijk}=\rho _{ijk}e^{-\delta _{i}\frac{%
 r_{i}^{4}}{t_{j}}}  \label{eq49}
 \end{equation}
 with $\delta _{i}=-\frac{\ln \left( P_{fa}\right) }{\alpha \left( \cos
 \left( \theta \right) \right) ^{2}}.$

 Since resolution cells are independent, the detection probability of the
 target $k$ in the direction $j$ is the sum of probabilities in the cells of
 this direction: 

 \begin{equation}
 P_{djk}=\sum_{i}^{N_{r}}P_{dijk}  \label{eq50}
 \end{equation}

 Lastly, we determine the probability $P_{dj}$ of detecting from one
 to several targets in a same direction. This probability is the union of
 previous probabilities for $k$ from one to $P.$ The Poincar\'{e} formula allows us to realize the calculation: 

 \begin{equation}
 \begin{array}{lll}
 P\left( \cup H_{i},i\in \left\{ 1,..,P\right\} \right) & = & 
 \sum_{i=1}^{n}P\left( H_{i}\right) - \\ 
 &  & \sum_{i,j=1,i\neq j}^{n}P\left( H_{i}\cap H_{j}\right) + \\ 
 &  & \sum_{i,j,l=1,i\neq j\neq l}^{n}P\left( H_{i}\cap H_{j}\cap H_{l}\right)
 \\ 
 &  & -...
 \end{array}
 \label{eq51}
 \end{equation}

 Our aim is to optimize this probability over the whole space.
 Unfortunately, its expression (\ref{eq51}) is not easily exploitable
 if we want to use the results described in the section above. It is
 the reason why we propose to realize a parametric modelling of this
 probability. The model we used is:
 
 \begin{equation}
 P_{dj}\simeq \exp \left( -\omega _{j}t_{j}^{-n_{j}}\right)  \label{eq52}
 \end{equation}

 where $\omega _{j}$ and $n_{j}$ are the modelling parameters in the
 direction $j$.  They are determined in order to minimize mean square
 error criterion. Under this formulation, the probability has the same
 properties as the one given by relation (\ref{eq4}); it is then
 possible to optimize it like using the framework of section
 \ref{capteur}, i.e. by a decomposition into an optimal number of
 elementary detections. Leading the same optimization process as in
 section \ref{capteur}, the following result can be shown. Let us call
 $\gamma_{sj}$ the unique solution of the equation :

 \begin{equation}
 \left( 1- \exp \left( - \gamma_{sj} \right) \right) \ln \left( 1- \exp \left( -\gamma_{sj} \right) \right) + n_j \gamma_{sj} \exp \left(\gamma_{sj}\right) =  0
 \end{equation}

 and $M_j$ the number of independant detections realized in the
 direction $j$ during an time $t_j$. If each elementary detection last
 $\frac{t_j}{M_j}$ then the detection probability from one to $N$
 targets in the direction $j$ is maximum when :

 \begin{equation}
 M_{j,opt}=\left( \frac{\gamma_{sj} t_j^{n_j}}{\omega_j} \right)^{\frac{1}{n_j}}.
 \end{equation}

 The elementary detection probability is then equal to $\exp(\gamma_{sj})$ and the overall detection probability to :
 
 \begin{equation}
 P_{dj}\left( t_{j} \right) =1-\exp \left( -\frac{t_{j}}{\tau _{j}}\right)
 \label{eq53}
 \end{equation}

 with :

 \begin{equation}
 \tau_j=-\left( \frac{\omega_j}{\gamma_{sj}} \right)^\frac{1}{n_j}\frac{1}{\ln \left( 1 - \exp \left( -\gamma_{sj} \right) \right)}.
 \end{equation}

 Using these last expressions for detection probabilities, the criterion to optimize is:

 \begin{equation}
 \mathcal{G}\left\{ 
 \begin{array}{ll}
 \text{maximize} & 
 \begin{array}{c}
 G=\sum\limits_{i=1}^{N_{d}}\varepsilon _{i}P_{dj}\left( t_{j}\right)
 \end{array}
 \\ 
 \text{under the constraints} & \left\{ 
 \begin{array}{l}
 \sum\limits_{i=1}^{N_d}t_{j}=T\quad \\ 
 t_{j}\geq 0\quad \forall j\in \left\{ 1,..,N_d \right\}
 \end{array}
 \right.
 \end{array}
 \right.  \label{eq14a}
 \end{equation}

 The resolution is similar as the one described in the deterministic context. It leads to the following results. Let $\lambda $ be the unique solution of the equation:
 
 \begin{equation}
 \sum_{i=1}^{N_d}\frac{\tau _{j}}{T}\left\lfloor \ln \left( \frac{T\varepsilon
 _{j}}{\tau _{j}\lambda }\right) \right\rfloor ^{+}-1=0  \label{eq22a}
 \end{equation}

 and $\mathcal{I}$ the suffix set defined by:
 
 \begin{equation}
 \mathcal{I}=\left\{ j \in \left\{ 1,..,N_d \right\} \left| \lambda < \frac{T\varepsilon _{j}}{\tau _{j}}\right. \right\}
 \label{eq23a}
 \end{equation}

 then the optimal temporal allocation is given by :
 
 \begin{equation}
 \left\{ 
 \begin{array}{lll}
 t_{j} & = & \tau _{j}\ln \left( \frac{T\varepsilon _{j}}{\tau _{j}\lambda } \right) \;if\;i\in \mathcal{I} \\ 
 & = & 0\;else
 \end{array}
 \right.  
 \label{eq24a}
 \end{equation}

 Furthermore, the optimal elementary detection number for each $t_{j}$, i.e for each direction, is equal to:

 \begin{equation}
 m_{j}= \left( \frac{\gamma_{sj} t^{n_j}}{\omega_j} \right)^\frac{1}{n_j}.
 \label{eq25a}
 \end{equation}

 \subsection{Simulations}
\begin{figure}
  \centering
  \resizebox{6cm}{!}{\includegraphics{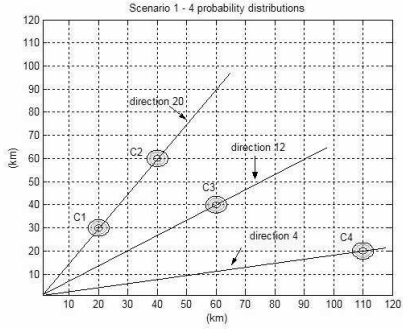}}
  \caption{\textit{Scenario 1: targets are located by distributions centers,straight lines represent the space directions}}
  \label{figsimul1}    
 \end{figure}

 We consider four targets located in an aerian space. An \textit{a
   priori} knowledge is available for each of them, by the way of
 density probabilities. These densities are defined in a cartesian
 frame by the mean of two dimensionnal gaussian laws. Figure
 \ref{figsimul1} illustrates this possible scenario which corresponds
 to the following situation:

 \begin{itemize}
 \item target 1 : the distribution is centered around the point
   $\left(20\,km,30\,km\right)$, the standard deviation on each
   coordinate is equal to $0.1\,km$.

 \item target 2 : the distribution is centered around the point
   $\left(40\,km,60\,km\right)$, the standard deviation on each
   coordinate is equal to $0.1\,km.$

 \item target 3 : the distribution is centered around the point
   $\left(60\,km,40\,km\right)$, the standard deviation on each
   coordinate is equal to $0.1\,km.$

 \item target 4 : the distribution is centered around the point
   $\left(110\,km,20\,km\right)$, the standard deviation on each
   coordinate is equal to $0.1\,km.$
 \end{itemize}

 According to our numerics values, space is divided into forty angular
 directions, $N_{d}=40$. The result of the optimization process leads
 to the allocation of table \ref{tab3} for $T=30\;ms$.

 \begin{table}[htbp] \centering

 \begin{tabular}{|c|c|c|c|c|c|c|c|}
 \hline
 {\small dir.} & {\small 1..3} & {\small 4} & {\small 5..11} & {\small 12} & 
 {\small 13..19} & {\small 20} & {\small 21..40} \\ \hline
 $\epsilon _{j}$ & {\small 1} & {\small 1} & {\small 1} & {\small 1} & 
 {\small 1} & {\small 1} & {\small 1} \\ \hline
 {\small t}$_{j}${\small \thinspace \thinspace \thinspace (ms)} & {\small 0}
 & {\small 0} & {\small 0} & {\small 23,47} & {\small 0} & {\small 6,53} & 
 {\small 0} \\ \hline
 {\small m}$_{j}$ & {\small 0} & {\small 0} & {\small 0} & {\small 1,31} & 
 {\small 0} & {\small 2,60} & {\small 0} \\ \hline
 {\small P}$_{dj}${\small \thinspace \thinspace \thinspace } & {\small 0} & 
 {\small 0} & {\small 0} & {\small 0,59} & {\small 0} & {\small 0,97} & 
 {\small 0} \\ \hline
 \end{tabular}

 \caption{\textit{Optimal temporal allocation without weights}\label{tab3}}%

 \end{table}

 As we can see only two directions are effectively considered in the
 temporal allocation. This is due to the allocation processus which
 realizes a global optimization of the detection probability. In fact,
 the fourth target is too distant from the others and from the sensor
 and it would spent too much time to detect it. This time would be
 allocated to the detriment of the other targets. It is therefore
 better not to observe it. Coefficients $\epsilon _{j}$ which appear in
 table \ref{tab3} are weights which introduce ponderations on the
 importance of the target.  They are all equal to one because the
 ponderation - or priority - notion was not taken into account
 initially. In table \ref{tab4} the result of the optimization process
 is written out when such taps are considered. As it can be seen in
 this table, the sensor spends more time in the important direction
 with respect to this taps. It results in an increasing in the
 detection probability.

 \begin{table}[htbp] \centering
 $
 \begin{tabular}{|c|c|c|c|c|c|c|c|}
 \hline
 $dir.$ & $1..3$ & $4$ & $5..11$ & $12$ & $13..19$ & $20$ & $21..40$ \\ \hline
 $\epsilon _{j}$ & $0$ & $0,07$ & $0$ & $0,18$ & $0$ & $0,74$ & $0$ \\ \hline
 $t_{j}\,\,\,(ms)$ & $0$ & $0$ & $0$ & $21,08$ & $0$ & $8,92$ & $0$ \\ \hline
 $P_{dj}\,\,\,$ & $0$ & $0$ & $0$ & $0,56$ & $0$ & $0,99$ & $0$ \\ \hline
 \end{tabular}
 $
 \caption{\textit{Optimal temporal allocation with weights}\label{tab4}}%

 \end{table}

 \section{Multisensor multitarget environment}

 \subsection{Introduction}
 In this section we suppose that $P$ radars are used to detect $N$
 targets, each sensor realizing a detection. The problem here is more
 complex than previously since we are looking for grouping radars in
 order to optimize the detection process. The questions to solve are
 therefore:

 \begin{itemize}
 \item which criterion do we want to optimize?
 \item how can we realize such a regrouping in a dynamical way? 
 \item how long a group a sensor must observe a potential target?
 \item how do we assign a group to a direction of observation?
 \end{itemize}

 This problem being tackled here in the deterministic context, the
 answer to the first question is therefore rather easy since it is the
 same as in the previous sections : the criterion to optimize is the
 sum of the detection probability :

 \begin{equation}
 \mathcal{C} = \sum_{k=1}^{N_t}p_k
 \label{equ:criterion3}
 \end{equation}

 with $N_t$ the target number and $p_k$ the detection probability of target $k$.

 In the following each group of sensors will be called a pseudo sensor.
 The detection of each pseudo sensor is derived from the fusion of the
 detection of each its sensors. We choose the fusion law OR which is
 usually used in the detection theory.

 The method proposed to answer to these questions is an heuristic based
 on the results of the previous sections. This heuristic is broken up
 into two phases :

 \begin{enumerate}
 \item The initial phase where first peudo sensors are constituted and
   where a first time allocation is made; this phase is based on the
   results of previous sections.
 \item The planification phase where the used of the sensors is
   planified over the time of analysis $T$ from the allocation realized
   during the initial phase.
 \end{enumerate}

 \subsection{The initial phase}
 The allocation process at initial time can be split into three steps:

 \begin{itemize}

 \item \textbf{Step 1}: Computation of the detection probabilities of
   the $P$\ sensors $K_i$ with $i \in \left\{ 1,...,P \right\}$.
   Knowing the ranges between the sensors and the targets and the
   observation duration $T$, it is possible, according to results of
   section \ref{hyp1}, to derive an optimal allocation over the
   duration $T$ for all the targets and the associated detection
   probabilities. This step allows the use of each sensor in an optimal
   way at initial time.

 \item \textbf{Step 2}: From the probabilities found at the step 1,
   compute the detection probabilities of each targets for all the
   possible pseudo-sensors. The aim of this step is to use the set of
   sensors in the best possible way at initial time. At this step the
   data fusion law is OR. To ensure that none of the sensors will be
   useless, the detection probability of each pseudo-sensor is computed
   with the shortest observation duration. For two sensors $A$ and $B$
   it comes to compute :
 
 \begin{equation}
 \begin{array}{lll}
 P\left( d_{A},t_{A},d_{B},t_{B}\right) & = & P_{d}\left( d_{A},\min \left(
 t_{A},t_{B}\right) \right)\\
 &  & +P_{d}\left( d_{B},\min \left( t_{A},t_{B}\right)\right) \\ 
 &  & -P_{d}\left( d_{A},\min \left( t_{A},t_{B}\right) \right)\\
 & &  P_{d}\left(d_{B},\min \left( t_{A},t_{B}\right) \right)
 \end{array}
 \end{equation}

 \item  \textbf{Step 3}: Determination of the allocation which maximize the criterion $\mathcal{C}$. 

 \end{itemize}

 Let us illustrate with an example the different steps of the method.

 \subsection{Example of initial allocation
 \label{exemplealloc}}

 Let us consider three sensors and three targets denoted by $C_{n},$
 $n\in \left\{ 1,2,3\right\}$. The table \ref{tab38} gives the
 distances between the sensors and the targets. For the sake of
 simplicity, angles $\theta$ are supposed to be null.

 \begin{table}[htbp] \centering
 \begin{tabular}{|c|c|c|c|}
 \hline
 $distance$ $sensor-target$ $\left( in\,\,km\right) $ & $K_{1}$ & $K_{2}$ & $%
 K_{3}$ \\ \hline
 $C_{1}$ & $45$ & $26$ & $52$ \\ \hline
 $C_{2}$ & $51$ & $45$ & $25$ \\ \hline
 $C_{3}$ & $50$ & $33$ & $41$ \\ \hline
 \end{tabular}
 \caption{\textit{Distances sensors-targets (in km)}\label{tab38}}
 \end{table}

 \subsubsection{Step 1: computation of the detection probabilities}

 Let us consider a duration $T=5\,ms$\ allocated to the detection
 phase. The optimal time allocation at initial time is written out in
 table \ref{tab39} and the corresponding detection probabilities in the
 table \ref{tab40}.

 \begin{table}[htbp] \centering
 \begin{tabular}{|c|c|c|c|}
 \hline
 $temporal$ \thinspace $allocation$ $\left( in\,\,s\right) $ & $K_{1}$ & $%
 K_{2}$ & $K_{3}$ \\ \hline
 $C_{1}$ & $2.5807$ & $1.1702$ & $0.9224$ \\ \hline
 $C_{2}$ & $1.0109$ & $1.8768$ & $1.1462$ \\ \hline
 $C_{3}$ & $1.4084$ & $1.9530$ & $2.9314$ \\ \hline
 \end{tabular}
 \caption{\textit{Optimal temporal allocation (in ms)}\label{tab39}}
 \end{table}

 \begin{table}[htbp] \centering
 \begin{tabular}{|c|c|c|c|}
 \hline
 $detection probabilities$ & $K_{1}$ & $K_{2}$ & $K_{3}$
 \\ \hline
 $C_{1}$ & $0.4814$ & $0.9309$ & $0.1233$ \\ \hline
 $C_{2}$ & $0.1444$ & $0.3797$ & $0.9532$ \\ \hline
 $C_{3}$ & $0.2095$ & $0.8206$ & $0.6612$ \\ \hline
 \end{tabular}
 \caption{\textit{Detection probabilities obtained from the data of the tables \ref{tab38} and \ref{tab39} }\label{tab40}}
 \end{table}

 \subsubsection{Step 2: pseudo-sensors and detection probabilities\label{pseudosensors}}

 Let us denote $P$ the number of sensors, then the number of
 pseudo-sensors is $S=2^{P}-1$. If $K_{1},$ $K_{2}$ and $K_{3}$ are the
 sensors, the pseudo-sensors are: $K_{1},$\ $K_{2},$\ $K_{3},$\
 $K_{1}-K_{2},$\ $K_{1}-K_{3},$\ $K_{2}-K_{3}$\ et $K_{1}-K_{2}-K_{3}.$
 The detection probabilities obtained by using the method described in
 the previous paragraph are given in the table \ref{tab41}.

 \begin{table}[htbp] \centering
 \begin{tabular}{|c|c|c|c|c|c|c|c|}
 \hline
   & {\small \textit{K}}$_{1}$ & {\small \textit{K}}$_{2}$ & {\small \textit{K}}$_{1}${\small -\textit{K}}$_{2}$ & {\small \textit{K}}$_{3}$ & {\small \textit{K}}$_{1}${\small -\textit{K}}$_{3}$ & {\small \textit{K}}$_{2}${\small -\textit{K}}$_{3}$ & {\small \textit{K}}$_{1}${\small -\textit{K}}$_{2}${\small -\textit{K%
 }}$_{3}$ \\ \hline
 {\small \textit{C}}$_{1}$ & {\small 0.481} & {\small 0.931} & {\small 0.949} & {\small 0.123} & {\small 0.307} & {\small 0.893} & {\small 0.916} \\ \hline
 {\small \textit{C}}$_{2}$ & {\small 0.144} & {\small 0.380} & {\small 0.338} & {\small 0.953} & {\small 0.943} & {\small 0.965} & {\small 0.956} \\ \hline
 {\small \textit{C}}$_{3}$ & {\small 0.210} & {\small 0.821} & {\small 0.771} & {\small 0.661} & {\small 0.530} & {\small 0.913} & {\small 0.864} \\ \hline
 \end{tabular}
 \caption{\textit{Detection probabilities associated to the pseudo-sensors}\label{tab41}}
 \end{table}

 \subsubsection{Step 3: determination of the optimal allocation}

 Using the results of section \ref{pseudosensors} it is easy to compute
 the value of the criterion (\ref{equ:criterion3}) for each possible
 allocation Target - Pseudo-sensor. The list below gives the results of
 the computation for a few pseudo-sensors. The allocation of a
 pseudo-sensor to a target is represented by an arrow.

 \begin{itemize}
 \item $K_{1}\rightarrow C_{1},$ $K_{2}\rightarrow C_{3}$ et $%
 K_{3}\rightarrow C_{2}$ : $0.4814+0.8206+0.9532=2.2552,$

 \item $K_{1}\rightarrow C_{3},$ $K_{2}\rightarrow C_{2}$ et $%
 K_{3}\rightarrow C_{1}$ : $\ 0.2095+0.3797+0.1233=0.7125,$

 \item $K_{2}\rightarrow C_{1},$ $K_{1}-K_{3}\rightarrow C_{3}$ : $\
 0.9309+0.5300=1.4609,$

 \item $K_{3}\rightarrow C_{1},$ $K_{1}-K_{2}\rightarrow C_{2}$ : $\
 0.1233+0.3384=0.4617,$
 \item $K_{1}-K_{2}-K_{3}\rightarrow C_{1}$ : $\ 0.9156,$
 \item $K_{1}-K_{2}-K_{3}\rightarrow C_{2}$ : $\ 0.9555,$
 \end{itemize}

 If we consider all the possible allocations, the maximum is obtained
 with the allocation which corresponds to the allocation of the sensor
 $K_{1}$ with the target $C_{1}$, the allocation of the sensor $K_{2}$
 with the target $C_{3}$ and the allocation of the sensor $K_{3}$ with
 the target $C_{2}$. This allocation, is such that all the targets are
 observed and all the sensors are used. It is interesting to remark
 that it is not the nearest sensor to a target which is used for its
 detection. This is because the context is not the optimization of the
 detection performances of each individual sensor but a context of
 global optimization.

 \subsection{Sensor planification over $T$}

 At this step, the initial allocation is realized. It is now necessary to
 build a planning of the use of the sensors from this initial
 allocation. If we analyze the initial allocation processus, we can see
 that this allocation is made for a given time interval resulting from
 the optimization process of the step 1 and from the limitation on the time
 of observation introduced in the step 2. This allocation is therefore
 not valid over the all time interval $T$. We propose the following
 planification of the sensor use over $T$ based on the results of the
 optimization process carried out during the step 1 of the initial
 allocation.

 \begin{itemize}
 \item rule 1 : the allocation sensor-target is called into question as
   soon as one of the durations of observation of the "active" couples
   is finished: the sensor or the sensors concerned must be allocated
   to another target. An active couple is an allocation of a sensor, or
   pseudo-sensor to a target.

 \item rule 2 : if a ponderation of the target has been achieved, the
   sensors are allocated to the target which have the highest
   ponderation weigths.

 \item rule 3 : if the are no ponderation, the sensors are allocated to
   the target which needs the lowest observation time different from
   zero. These times are those computed at the the step 1 of the
   initial allocation. Thus, by given priority to the short durations,
   the detection performances of the observed targets will be optimized
   in the case where some operational constraints abort the detection
   process.

 \end{itemize}

 \subsection{Exemple of planification}
 
 We consider in this example the situation used in the section
 \ref{exemplealloc} with all the ponderation taps equal to one. The
 allocation sensors-targets has been determined in the section
 \ref{exemplealloc}. The duration during which this allocation is
 effective is determined by rule 1. It corresponds to the minimum
 observation duration of the targets by the selected sensors.
 Considering the results of table \ref{tab39}, the allocation is then
 called into question at the end of the time $t_{obs}=1.1462\,ms$.
 During this duration, the sensor $K_{3}$ has observed the target
 $C_{2}$ in an optimal way, in the sense of the detection performances,
 and it then needs to be directed towards another target. The other
 targets have also been observed during $t_{obs}$. The observation of
 $C_{2}$ by $K_{3}$ being finished the allocation given by table
 \ref{tab39} must be modified to make appear that $K_3$ will no more
 observe $C_2$ and that for the other affectation targets have still
 being observed during $t_{obs}$. The result is given is written out in
 the table \ref{tab42}.
 
 \begin{table}[htbp] \centering
 \begin{tabular}{|c|c|c|c|}
 \hline
 $temporal$ \thinspace $allocation$ $\left( in\,\,km\right) $ & $K_{1}$ & $
 K_{2}$ & $K_{3}$ \\ \hline
 $C_{1}$ & $\mathbf{1.4345}$ & $1.1702$ & $0.9224$ \\ \hline
 $C_{2}$ & $1.0109$ & $1.8768$ & $\mathbf{0}$ \\ \hline
 $C_{3}$ & $1.4084$ & $\mathbf{0.8068}$ & $2.9314$ \\ \hline
 \end{tabular}
 \caption{\textit{Temporal allocation after an observation duration $t_{obs}$}\label{tab42}}
 \end{table}

\begin{figure}[ht]
 \begin{center}
   \includegraphics[width=6cm]{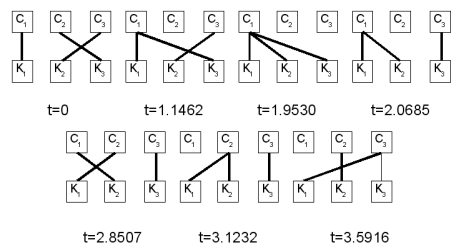}
   \caption{\textit{Representation of the planning of the utilization of the sensors during the duration }$T$}
   \label{figalloc}
 \end{center}
 \end{figure}

 Using the rule 3, the sensor $K_{3}$ is now oriented towards the
 target, which need the lowest observation duration in order to reach
 optimal detection performances. Here, it is the target $C_{1}$. The
 sensor $K_{2}$ remains affected to the observation of the target
 $C_{3}$, then the pseudo-sensor $K_{1}-K_{3}$ is used for the
 observation of the target $C_{1}$. This algorithm is iterated till the
 time $t=T$ is reached. The resulting allocation is represented in the
 figure \ref{figalloc}. The target $C_{1}$ will have been observed
 during $3.1232ms$, the target $C_{2}$ during $3.5655ms$ and the target
 $C_{3}$ during $4.8845ms$. Their detection probabilities are
 respectively $0.9686,$ $0.9751$ and $0.9520$. The sum of these
 probabilities is $2.8957$. It it is greater than $2.7075$ which is the
 sum which would have been obtained if the allocation established
 initially had been used during all the duration $T$.

 \section{Conclusion}
 This paper presents methods to manage the time allocation of radars
 over a set of targets. In a first part a method to optimize the
 detection process of targets is proposed. It is based on the modelling
 of the detection probability of a target. This firt result is then
 used to propose optimal time allocations in the monosensor multitarget
 case. Two operational contexts are considered : a deterministic
 context where the position of the target are known and a probabilistic
 context where the knowledge of the position of the target is
 represented by probability density functions. We showed that the
 probabilistic context can be solved using the results of the
 deterministic one. These results have then been used to propose an
 heuristic for the planification of a set radar the mission of which is
 to detect targets : we are then in the multisensor multitarget case.
 The planification has been proposed in the deterministic context and
 still need to be generalized to the probabilistic context.

 \section{Acknowledgement}
 The authors would like to thank Michel Prenat from Thales Optronics for its important contribution to this work.

\bibliographystyle{IEEEtran}
\bibliography{m2vbib,MultiSensor}






\end{document}